\documentclass[10pt]{article}
\usepackage{amssymb}

\usepackage{enumerate}
\usepackage{graphicx}
\usepackage{amsthm}
\usepackage{color}
\usepackage{epstopdf}

\def\R{{\rm I\! R}}
\def\N{{\rm I\! N}}

\newtheorem{theorem}{Theorem}
\newtheorem{corollary}{Corollary}

\newtheorem{lemma}{Lemma}

\newtheorem{remark}{Remark}

\title{Limit cycle's uniqueness for second order O.D.E.'s polynomial in $\dot x$}

\author{M. Sabatini  
\footnote{Dip. di Matematica, Univ. di Trento, I-38050 Povo, (TN) - Italy.
Email: marco.sabatini@unitn.it,
Phone: ++39(0461)881670, Fax: ++39(0461)881624 }
}
\date{}
\begin{document}
\maketitle
\begin{abstract}
We prove a uniqueness result for limit cycles of the second order ODE $\ddot x + \sum_{j=1}^{J}f_{j}(x)\dot x^{j} + g(x) = 0$.  We extend a uniqueness result proved in \cite{CRV}. The main tool applied is an extension of Massera theorem proved in \cite{GS}.

{\bf Keywords}: Uniqueness, limit cycle, second order ODE's, Massera theorem, Conti-Filippov transformation.
\end{abstract}

\section{Introduction}
 
In this paper we are concerned with planar differential systems of the form
\begin{equation}\label{sysypol} 
\dot x =  y, \qquad \dot y = - g(x) -  \sum_{j=1}^{J}f_{j}(x)y^{j},
\end{equation} 
equivalent to the second order differential equations of the form
\begin{equation}\label{equaypol} 
\ddot x + \sum_{j=1}^{J}f_{j}(x)\dot x^{j} + g(x) = 0.
\end{equation} 
Several mathematical models of physics, economics, biology are governed by second order differential equations (\cite{LC}, \cite{SC}, \cite{V}). Other models can be reduced to systems of the type (\ref{sysypol}) by means of suitable transformations. 
The asymptotic behaviour of their solutions is one of the main objects of study. In this perspective, the existence of special solutions as stationary ones, or isolated cycles, is of primary interest. This is particularly true if such solutions attract (repel) neighbouring ones, so that the system's dynamics is dominated by that of the attracting equilibria or cycles. In the special case of  an isolated cycle attracting all the other solutions but equilibria, the description of the system's dynamics becomes quite simple, since the asymptotic behaviour of  all solutions but the equilibrium one is just that of the limit cycle. 
Uniqueness theorems for limit cycles have been  extensively studied (see \cite{CRV}, \cite{XZ1}, \cite{XZ2}, \cite{Ci} for recent results and extensive bibliographies, \cite{ZDHD}, chapter IV, section 4). In general, studying the number and location of  limit cycles is a non-trivial problem, as shown by the resistance of   Hilbert XVI problem. Such a problem has been recently re-proposed as a main research problem (see \cite{Sm}, problem 13). A strictly related subject is that of hyperbolicity. A $T$-periodic cycle $\gamma(t)$ of a differential system
\begin{equation}\label{sysPQ} 
\dot x =  P(x,y), \qquad \dot y =  Q(x,y),
\end{equation} 
is said to be {\it hyperbolic} if
\begin{equation}\label{hyper} 
\int_0^T {\rm div } (\gamma(t)) dt \neq 0,
\end{equation} 
where ${\rm div } = {\partial P \over \partial x} +  {\partial Q \over \partial y}$ is the divergence of $(\ref{sysPQ})$. Hyperbolicity plays a main role in perturbation problems, since smooth perturbations of hyperbolic cycles do not allow multiple bifurcations. An attractive cycle is not necessarily hyperbolic. 

Most of the uniqueness results proved for planar systems are concerned with the classical Li\'enard system  and its generalizations, such as
\begin{equation}\label{sysGG}
\dot x = \xi(x)\bigg[ \varphi(y) - F(x)  \bigg], \qquad \dot y = -\zeta(y)g(x), \quad \xi(x) \neq 0, \quad \zeta(y) \neq 0.
\end{equation}
Such a class of systems also contain Lotka-Volterra systems and systems equivalent to Rayleigh equation as special cases \cite{LC}.  Such systems are characterized by the presence, both in $\dot x$ and $\dot y$, of a single mixed term obtained as the product of single-variable functions, resp.   $ \xi(x) \varphi(y) $ and $\zeta(y)g(x)$. Moreover, such systems can be easily transformed into systems without mixed terms, by applying the transformation
$$
X(x) = \int_0^x \frac{1}{\xi(s)} ds, \qquad Y(y) = \int_0^y \frac{1}{\zeta(s)} ds. 
$$
The transformed system has the form
$$
\dot X = \tilde{ \varphi}(Y) - \tilde {F}(X) , \qquad \dot y =  \tilde {g}(X),
$$
for suitable functions $ \tilde{ \varphi}(Y)$, $\tilde {F}(X)$, $ \tilde {g}(X)$.  \\
\indent In order to study systems with several distinct mixed terms, a different approach is required. 
 Some recent results (\cite{CRV}, \cite{Ci})  are concerned with the following systems, 
\begin{equation}\label{sysCRV} 
\dot x = y, \qquad \dot y = - x -  y \sum_{k=0}^{N}f_{2k+1}(x)y^{2k},
\end{equation} 
equivalent to the equations
\begin{equation}\label{equaCRV} 
\ddot x +  \sum_{k=0}^{N}f_{2k+1}(x){\dot x}^{2k+1} + x = 0.
\end{equation} 
Such systems cannot be reduced to the form (\ref{sysGG}) by the above transformation. In order to prove the limit cycle uniqueness, the authors extend a classical result by Massera about the uniqueness and global  (except for equilibria) attractiveness of Li\'enard limit cycles \cite{M}. Such a result comes from two main properties. The first one is that every cycle be star-shaped, fact proved in a new way in \cite{CRV}. The second one is that the vector field rotates clockwise along rays (half-lines having extreme at the axes' origin $O$).  
In fact, in general such properties are not sufficient to prove the limit cycle's uniqueness, as the following example shows, 
\begin{equation}\label{sysduecicli} 
\left\{
\begin{array}{rl}
\dot x   = & y \Big(x^2+y^2-(x^2+y^2)^2 \Big)  + x \Big(1-3 (x^2+y^2)+(x^2+y^2)^2 \Big) \\ 
 \dot y  = & -x \Big(x^2+y^2-(x^2+y^2)^2 \Big) + y \Big(1-3 (x^2+y^2)+(x^2+y^2)^2 \Big)
\end{array}
\right.
\end{equation} 
Such a system has two star-shaped limit cycles coinciding with the circles $x^2+y^2 =\frac{3-\sqrt {5}}{2}$ and  $x^2+y^2 =\frac{3+\sqrt {5}}{2}$.  The internal one is an attractor, the external one is a repellor. The vector field rotates clockwise along every ray (see figure 1). 
\begin{figure}[h!]
  \caption{The system (\ref{sysduecicli}) has two limit cycles.}
  \centering
    \includegraphics[width=0.5\textwidth]{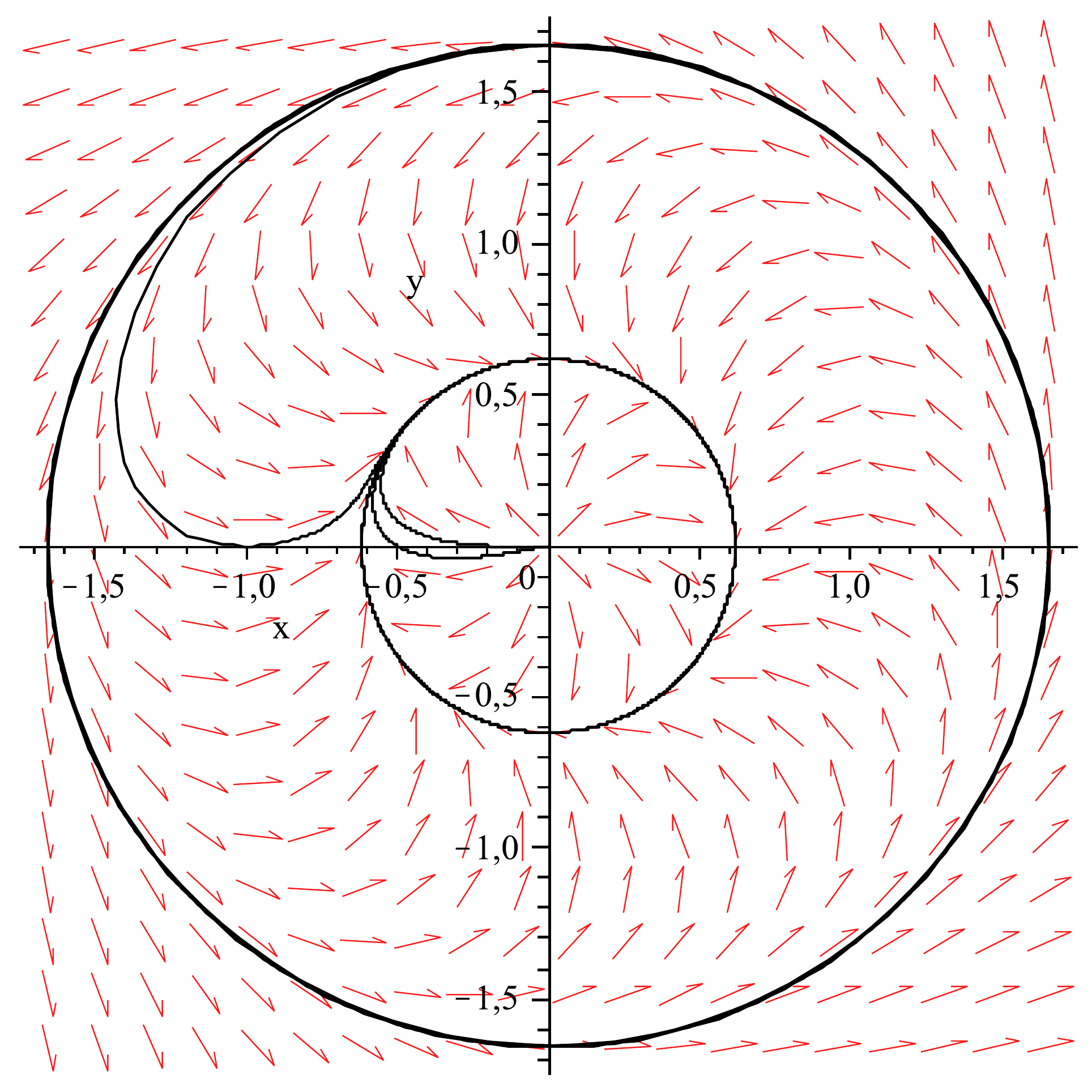}
\end{figure}

\indent  An even more pathological example (from the point of view of Massera-like theorems) is the system
\begin{equation}\label{systrigo} 
\dot x  = y \cos (x^2+y^2) - x \sin (x^2+y^2), \quad \dot y  = -x \cos (x^2+y^2) - y \sin (x^2+y^2).
\end{equation} 
Such a system satisfies both properties, and  has infinitely many cycles coinciding with the circles $\{ x^2+y^2 = k \pi : k \in \N \}$. Such cycles are alternatively repelling  and attracting, and rotate alternatively  counter-clockwise and clockwise. A cycle's attraction (repulsion) region is the annular region bounded by the two adjacent limit cycles, except for the innermost one, whose attraction region has the origin in its boundary. \\
\indent  As a consequence, an additional condition is required in order to give a complete proof of a Massera-like theorem. Corollary 6 of \cite{GS} provides a natural additional condition, asking for the angular velocity not to vanish in the whole plane. 
The same result can be proved if the angular velocity does not vanish in a suitable subregion of the plane, as in \cite{S}. Corollary 6 is a particular case of   theorem 2 in \cite{GS}, a rather general extension to Massera theorem. In such a theorem an auxiliary function $\nu$ is used in order to study the cycles' hyperbolicity. The main property of $\nu$ is that its integral on a cycle coincides with that of the vector field's divergence. On the other hand, $\nu$ has an advantage over the divergence, since it can be everywhere negative in presence of a repelling critical point and an attracting limit cycle, while in such a situation the divergence has to change sign. This helps in proving the limit cycle's uniqueness, as in \cite{S}. Additionally, theorem 2 in \cite{GS} allows to prove the limit cycle's hyperbolicity, which is not  a consequence of Massera-like theorems.

In this paper we extend in two ways  the uniqueness result of \cite{CRV}. First, we consider systems with even degree terms $f_{2k}(x)y^{2k}$, then we allow  $g(x)$ to be non-linear, provided $xg(x) \neq 0$ for $x \neq 0$ in some interval.  Even degree terms can be treated by considering the function $\sum_{j=1}^{J}f_{j}(x)y^{j-1}$ as the sum of $y$-trinomials, each satisfying suitable conditions. On the other hand, considering a non-linear $g(x)$ allows to work in  different regions for the same system. Since a $x$-translation transforms a system of the form (\ref{sysypol}) into a system of the same form, one just has to translate  a critical point to the origin, proving the uniqueness, in a suitable strip $(a,b) \times \R$, of limit cycles surrounding such a point.    \\
\indent   If the system has the form (\ref{sysCRV}), with $g(x) = x$,  our result extends that one obtained in \cite{CRV}, replacing the sign and monotonicity conditions on $f_{2k+1}(x)$ with a monotonicity condition on $x^{2k}f_{2k+1}(x)$. This allows to apply our theorem to some polynomial coefficients, as $f_{3}(x) = x^4-x^2+1$, which do not satysfy the hypotheses in \cite{CRV}.   \\
\indent In order to extend the results of  \cite{CRV}, we follow a different approach w. r. to that one developed in \cite{S}. We cannot apply the theorem 1 of \cite{S} to \ref{sysCRV}, since $\phi(x,y) =  \sum_{k=0}^{N}f_{2k+1}(x)y^{2k}$ does not satisfy the strict star-shapedness condition required by such a theorem. In order to overcome such an obstacle we introduce a non-invariance property similar to that one appearing in LaSalle Invariance Principle. In fact, the passage from $x\phi_x + y\phi_y > 0 $ to $x\phi_x + y\phi_y \geq 0 $ is analogous to the passage from the condition $\dot V < 0$ to the condition $\dot V \leq 0$ in studying asymptotic stability problems.
Moreover, rather then just proving the cycle's star-shapedness, we prove that, assuming by absurd the existence of  two concentric limit cycles, then there exists an annular region containing both ones, where the angular velocity does not vanish. This allows to apply theorem 2 in \cite{GS} in order to get the limit cycle's uniqueness and its hyperbolicity. \\
\indent    Our paper is organized as follows. In section 1 we prove the main theorem about uniqueness and hyperbolicity  for systems with a linear $g(x)$. Then we apply such a theorem to some special cases. \\
\indent   In section 2 we apply Conti-Filippov transformation in order to reduce a system of the form (\ref{sysypol}) with a non-linear $g(x)$ to another system of the same type, but with a linear $g(x)$. In this section we derive the uniqueness condition in strips.  \\
\indent  In both sections we also consider the question of limit cycles existence, giving sufficient conditions for the solutions' boundedness that allow to apply Poincar\'e-Bendixson theorem.

\section{g(x) linear}

Since our approach is based on the cycle's star-shapedness, we restrict to a star-shaped subset $\Omega \subset \R^2$. In this context, we think of $\Omega$ as a strip $(a,b) \times \R$, with $a < 0 < b$, but what follows holds for arbitrary star-shaped subsets. We denote partial derivatives by subscripts, i. e. $\phi_x$ is the derivative of $\phi$ w. r. to $x$, etc..
We say that a function $\phi \in C^1( \Omega, \R) $ is  {\it star-shaped} if $(x,y) \cdot \nabla \phi = x \phi_x + y \phi_y$ does not change sign. We say that $\phi$ is {\it strictly star-shaped} if  $(x,y) \cdot \nabla \phi \neq 0$, except at the origin $O=(0,0)$. We call {\it ray} a half-line having origin at the point $(0,0)$.  
 We denote by $\gamma(t,x,y)$ the unique  solution of the system (\ref{sysypol}) such that $\gamma(0,x,y) = (x,y)$.  For definitions related to dynamical systems, we refer to \cite{BS}. 
 
 Throughout all of this paper we assume  the system (\ref{sysypol}) to satisfy some hypotheses ensuring existence and uniqueness of solutions, and continuous dependence on initial data. This occurs if, for instance, one has
 
 $\bullet )$ \qquad $f_{j}  $, $j= 1 , \dots, J$, continuous on their domains;
 
  $\bullet )$ \qquad $g$ lipschitzian on its domain.

In this section we are concerned with the system (\ref{sysypol}), assuming  $g(x)=kx$. 
Without loss of generality, possibly performing a time rescaling, we may restrict to the case $k=1$. In this case, we may consider (\ref{sysypol}) as a special case of a more general class of systems,
\begin{equation}\label{sysphi} 
\dot x = y, \qquad \dot y = - x - y\phi(x,y).
\end{equation} 
We first consider a sufficient condition for limit cycle's uniqueness. 
It has the same form as system (9) in \cite{S}, but here the function $\phi(x,y)$ is not strictly star-shaped. This will force us to modify the proof in the part related to the angular velocity, and add an hypothesis in the main  theorem. As in \cite{S}, we set
$$
A(x,y) = y \dot x - x \dot y = y^2 + x^ 2 +xy \phi(x,y).
$$
The sign of $A(x,y)$ is opposite to that of the angular speed of the solutions of (\ref{sysphi}).  
In general, $A(x,y)$ changes sign along the solutions of (\ref{sysphi}), even in the simplest example of nonlinear Li\'enard system with a limit cycle
\begin{equation}\label{sysVDP} 
\dot x = y, \qquad \dot y = - x - y(x^2-1).
\end{equation}
In fact, the orbits of (\ref{sysVDP}) approaching the limit cycle from the second and fourth orthants change angular velocity when they get closer to the limit cycle.  
 
Our uniqueness result comes from theorem 2 of \cite{GS}, in the form of corollary 6. In the following lemma the 
we prove that if two limit cycles $\mu_1$ and $\mu_2$  exist, then in the closed annular region $D_{12}$ bounded by $\mu_1$ and $\mu_2$ the function $A$ does not vanish. Next lemma's proof is a modification of part of the proof of theorem 1 in \cite{S}.  We emphasize that we consider {\it open} orthants, i. e. orthants without semi-axes.

\begin{lemma}\label{lemma1} Let $\phi\in C^1( \Omega, \R^2)$, with  $x\phi_x + y\phi_y \geq 0$ in $\Omega$. If the system (\ref{sysphi}) has two distinct limit cycles  $\mu_1$ and $\mu_2$, then $A(x,y) >0$ on  $D_{12}$.
\end{lemma}
{\it Proof.}
\indent The system (\ref{sysphi})  has one critical point, hence  $\mu_1$ and $\mu_2$ are concentric.  Let $\mu_1$ be the inner one, $\mu_2$ be the outer one. The radial derivative $A_r$ of $A$ is given by
$$
A_r = \frac {x A_x + y A_y}r  = \frac 1r\bigg(2 A + xy (x \phi_x + y\phi_y) \bigg) = \frac 1r\bigg( 2 A + xyr \phi_r \bigg).
$$
Hence the function $A$ satisfies the differential equation $rA_r =  2 A + xy (x \phi_x + y\phi_y) $, whose right-hand side  sign determines the monotonicity of $A$ along rays. If $x^*y^* > 0$ and $A(x^*,y^*) > 0$, then $A_r >0$ at $(x^*,y^*)$ and at every point $(rx^*,ry^*)$ with $r > 1$, hence $A$ is strictly increasing on the half-line $\{ (rx^*,ry^*) : r > 1\}$. On the other hand, if $x^*y^* < 0$ and $A(x^*,y^*) < 0$, then $A_r < 0$ at $(x^*,y^*)$ and at every point $(rx^*,ry^*)$ with $r > 1$, hence $A$ is strictly decreasing on the half-line $\{ (rx^*,ry^*) : r > 1\}$.   \\
\indent If $A(x^*,y^*) = 0$ and $x^*y^* >0 $, then either $A(rx^*,ry^*) = 0$ for all $r > 1$, or at some point of the half-line $\{ (rx^*,ry^*) : r > 1\}$ the radial derivative $A_r$ becomes positive, so that $A(rx^*,ry^*) > 0$ for all $r > \overline{r}$, for some $\overline{r} > 1$.  \\
\indent Similarly, if $A(x^*,y^*) = 0$ and $x^*y^* < 0 $, then either $A(rx^*,ry^*) = 0$ for all $r > 1$, or at some point of the half-line $\{ (rx^*,ry^*) : r > 1\}$ the radial derivative $A_r$ becomes negative, so that $A(rx^*,ry^*) < 0$ for all $r > \overline{r}$, for some $\overline{r} > 1$. \\
 \indent We prove  that, for every orbit $\gamma$ contained in $D_{12}$, $A(\gamma(t)) > 0$. Every orbit in $D_{12}$ meets every semi-axis, otherwise its positive limit set would  contain a critical point different from $O$. On every semi-axis one has $A(x,y) > 0$. Assume first, by absurd, $A(\gamma(t))$ to change sign. Then there exist $t_1 < t_2$ such that $A(\gamma(t_1)) > 0$, $A(\gamma(t_2)) < 0$, and and $\gamma(t_i)$, $i = 1,2$ are on the same ray. Assume $\gamma(t_i)$, $i = 1,2$ to be in the first orthant. Two cases can occur: either  $|\gamma(t_1)| < |\gamma(t_2)|$ or  $|\gamma(t_1)| > |\gamma(t_2)|$. The former, $|\gamma(t_1)| < |\gamma(t_2)|$, contradicts the fact that $A$ is radially increasing in the first orthant, hence one has $|\gamma(t_1)| > |\gamma(t_2)|$. The orbit $\gamma$ crosses the segment $\Sigma = \{r \gamma(t_1), 0 < r <1\}$ at $\gamma(t_2)$, going towards the positive $y$-semi-axis. Let $G$ be the sub-region of the first orthant bounded by the positive $y$-semi-axis, the ray $\{r \gamma(t_1), r > 0\}$ and the portions of $\mu_1$, $\mu_2$ meeting the $y$-axis and such a ray.  The orbit $\gamma$ cannot remain in $G$, since in that case $G$ would contain a critical point different from $O$. Also, $\gamma$ cannot leave $G$ crossing the positive $y$-semi-axis, because $A(x,y) > 0$ on such an axis. Hence $\gamma$ leaves $G$ passing again through the segment $\Sigma$. That implies the existence of $t_3 > t_2$, such that $\gamma(t_3)$ lies on the ray  $\{ r \gamma(t_1), 0 < r \}$. Again, one cannot have $|\gamma(t_3)| < |\gamma(t_2)|$, since $A(\gamma(t_3)) > 0$ implies $A$ increasing on the half-line $r \gamma(t_3), r > 1$, hence one has $|\gamma(t_3)| > |\gamma(t_2)|$. Also, one cannot have $|\gamma(t_3)| < |\gamma(t_1)|$, otherwise $\gamma$ would enter a positively invariant region, bounded by the curve $\gamma(t)$, for $t_1 \leq t \leq t_3$,  and by  the segment with extrema $\gamma(t_1)$, $\gamma(t_3)$, hence there would exist a critical point different from $O$. As a consequence, one has $|\gamma(t_3)| > |\gamma(t_1)|$.
Since $A(x,y) >0$ on the segment joining $\gamma(t_1)$ and $\gamma(t_3)$, such a segment, with the portion of orbit joining $\gamma(t_1)$ and $\gamma(t_3)$ bounds a region which is negatively invariant for (\ref{sysphi}), hence contains a critical point different from $O$, contradiction.  \\
\indent   This argument may be adapted to treat also the case of a ray in the second orthant, replacing the positive $y$-semi-axis with the positive $x$-semi-axis, and reversing the relative positions of the points $\gamma(t_1)$, $\gamma(t_2)$, $\gamma(t_3)$. In the second orthant one uses the fact that if $A(\gamma(t^*)) < 0$, then on the half-line $r A(\gamma(t^*)), r > 1$ the function $A$ is strictly decreasing, since $xy (x \phi_x + y\phi_y) \leq 0$.  In the third and fourth  orthants one repeats the arguments of the first and second orthants, respectively.  \\
\indent Finally, assume that $A(\gamma(t)) = 0$ at some point $(x^*,y^*) \in D_{12}$. First, we assume $(x^*,y^*)$ to be in the first orthant.  \\
\indent  One cannot have $A(x^*,y^*) = 0$ on all of the segment $(rx^*,ry^*), 0 < r < 1$, since in such a case the vector field would be radial, and  $\gamma$ would contain the whole segment, contradicting the fact that $\gamma \subset D_{12}$. Also, one cannot have $A(x^*,y^*) > 0$ at any point $(x^+,y^+)$ of the segment $(rx^*,ry^*), 0 < r < 1$, since in that case $A$ would be increasing on the half-line $(rx^+,ry^+), r > 1$, contradicting $A(x^*,y^*) = 0$. Hence there exists a point $(x^-,y^-)$ of the segment $(rx^*,ry^*), 0 < r < 1$ such that $A(x^-,y^-) < 0$. Moreover, for all $0<r<1$ one has  $A(rx^-,ry^-) < 0$. 
Working in the same way one proves that the segment $(rx^*,ry^*), 0 < r < 1$ is the disjoint union  of an open sub-segment $\Sigma^-$ where $A(x,y) < 0$, and a closed sub-segment  $\Sigma^0$ where $A(x,y) = 0$.
Then the "triangle" $T$ bounded by the $y$-semi-axis, the ray $(rx^*,ry^*), r> 0$, and the arc of $\mu_1$, is positively invariant. This produces a contradiction, due the fact that orbits starting in $T$, close to $\mu_1$, by the continuous dependence on initial data have to remain close to $\mu_1$, hence have to get out of $T$. \\
\indent  Now we assume $(x^*,y^*)$ to be in the second orthant and reverse the above argument: one cannot have $A(x^*,y^*) = 0$ on all of the half-line $(rx^*,ry^*), r > 1$, since in such a case $\gamma$ would contain the whole half-line, contradicting the fact that $\gamma \subset D_{12}$.  Let $r^* > 0$ be the maximum positive $r$ such that  $A(rx^*,ry^*) = 0$. Then for all $r > r^*$, one has $A(rx^*,ry^*) \neq 0$, hence, by what said at the beginning of this proof,  $A(rx^*,ry^*) < 0$ for $r > r^*$. Hence, the "unbounded triangle" $T$ having as boundary part of the positive $x$-semi-axis, the half-line $(rx^*,ry^*), r > r^*$, and the arc of $\mu_2$ connecting such half-lines, is positively invariant. As above, this produces a contradiction, due the fact that orbits starting in $T$, close to $\mu_2$, by the continuous dependence on initial data have to remain close to $\mu_2$, hence have to get out of $T$. \\
\indent In the other two orthants one can repeat the arguments of the first two ones, completing the proof. 
\hfill$\clubsuit$

Next theorem is as well a modification of theorem 1 in \cite{S}. In order to cope with the weaker hypothesis, we introduce a new condition on the set $x\phi_x +y\phi_y =0$.

\begin{theorem}\label{teorema} Let $\phi\in( \Omega, \R^2)$ be a star-shaped function, such that the set $x\phi_x +y\phi_y =0$ does not contain any non-trivial positive semi-orbit of  (\ref{sysphi}). Then (\ref{sysphi}) has at most one limit cycle, which is hyperbolic.
\end{theorem}
{\it Proof.}
Let us assume, by absurd, (\ref{sysphi}) to have two distinct limit cycles $\mu_1$ and $\mu_2$. Since the system has only one critical point, $\mu_1$ and $\mu_2$ have to be concentric. Assume $\mu_1$ to be the inner one. Let us restrict to the closed annular region $D_{12}$ bounded by $\mu_1$ and $\mu_2$. By the lemma \ref{lemma1}, one has $A(x,y) > 0$ on $D_{12}$. Let us consider the new system obtained by dividing the vector field of (\ref{sysphi}) by $A(x,y)$, as in  corollary 6 in \cite{GS}. In order to apply such a corollary, one has to  compute the expression
$$
\nu = \frac{P\left( xQ_x + yQ_y \right) - Q \left( xP_x + yP_y\right) }{y P - x Q} 
$$
where $P$ and $Q$ are the components of the considered vector field. Since for system (\ref{sysphi}) one has $y P - x Q = A$ , one can write
$$
\nu A=
y \left(-x - xy \phi_x - y\phi - y ^2 \phi_y \right) - \left(  -x - y \phi(x,y) \right) y  =
$$
$$
-y^2 \left( x \phi_x + y  \phi_y \right) \leq 0.
$$The function $\nu$ vanishes for $y=0$ and for $x \phi_x + y  \phi_y= 0$. The set $y=0$ is transversal to both cycles, hence both cycles have two points on $y=0$. Moreover, since by hypothesis the set $x\phi_x +y\phi_y =0$ does not contain any positive semi-orbit of  (\ref{sysphi}), every cycle contains at least a point such that $x\phi_x +y\phi_y > 0$. By continuity, $x\phi_x +y\phi_y >0$ in a neighbourhood of such a point, so that $x\phi_x +y\phi_y >0$ on an arc of $\mu_i$, $i = 1 , 2$.  \\
\noindent Then, for both cycles one has:
$$
\int_{0}^{T_i} \nu(\mu_i(t) )dt < 0, \qquad i=1,2,
$$
where $T_i$ is the period of $\mu_i$, $i=1,2$. Hence both cycles, by theorem 1 in \cite{GS}, are attractive. Let $A_1$ be the region of attraction of $\mu_1$.  $A_1$ is bounded, because it is enclosed by $\mu_2$, which is not attracted to $\mu_1$. The external component of   $A_1$'s boundary is itself a cycle $\mu_3$, because (\ref{sysphi}) has just one critical point at $O$. We can apply to $\mu_3$ the above argument about the sets $y=0$ and $x\phi_x +y\phi_y =0$, concluding that
$$
\int_0^{T_3} \nu(\mu_3(t) )dt < 0.
$$
Hence $\mu_3$ is attractive, too. This contradicts the fact that the solutions of (\ref{sysphi}) starting from its inner side are attracted to $\mu_1$. Hence the system  (\ref{sysphi}) can have at most a single limit cycle. Its hyperbolicity comes from the equality (see \cite{GS})
$$
\int_0^T {\rm div } (\mu(t)) dt = \int_0^{T} \nu(\mu(t) )dt < 0.
$$
\hfill $\clubsuit$

Now we consider some classes of $y$-polynomial systems satisying the hypotheses of theorem \ref{teorema}. Let us consider the following conditions related to a $y$-trinomial $\kappa(x) y^{2h+2r} + \tau(x) y^{h+2r} + \eta(x)y^{2r}$:

${\bf (T^+)}$ \quad for all $x \in (a,b)$, $x \neq 0$, one has $(x \tau' + (h+2r)\tau)^2 - 4( x\eta' + 2r \eta)(x \kappa' + (2h+2r) \kappa)) \leq 0$, and $x \kappa' + (2h+2r) \kappa \geq  0$.

\begin{remark}Since the above inequality implies a sign condition on the $y^h$-trinomial $(x\kappa' + 2(h+r) \kappa )y^{2h} + (x \tau' + (h+2r) \tau ) y^{h} + (x\eta' + 2r\eta) $ (see next corollary), it is equivalent to require $x \kappa' + (2h+2r) \kappa \geq  0$ or  $x\eta' + 2r \eta \geq 0$.
\end{remark}

${\bf (Seq)}$ \quad   there exists a sequence $x_m$ converging to $0$, such that for every  $m$ there exists $1 \leq j(m) \leq N$ satisying $x_mf_{j(m)}'(x_m) + (j(m) - 1)f_{j(m)}(x_m) \neq 0$.

\begin{corollary}\label{corollarioT} Assume   $\sum_{j=1}^{J}f_{j}(x)y^{j-1}$ to be the sum of $y$-trinomials satisfying the conditions $(T^+)$ and $(Seq)$.   Then  the system (\ref{sysypol}) has at most one limit cycle, which is hyperbolic.
\end{corollary}
{\it Proof.}
The system (\ref{sysypol}) is a special case of the system  (\ref{sysphi}), with 
$$
\phi(x,y) = \sum_{j=1}^{J}f_{j}(x)y^{j-1}. 
$$
Computing $\nu$, one has
$$
\nu = -y^2 \left(x \phi_x + y  \phi_y  \right)   = - y^2  \sum_{j=1}^{J} \big(xf_{j}'(x) + (j-1) f_{j}(x) \big) y^{j-1}.
$$
By hypothesis, there exist  $y$-polynomials $P_l(x,y) = \kappa_l(x) y^{2h(l)+2r(l)} + \tau_l(x) y^{h(l)+2r(l)} + \eta_l(x)y^{2r(l)}$ such that
$$
\phi(x,y) = \sum_{j=1}^{J}f_{j}(x)y^{j-1} = \sum_{l=1}^{L} P_l(x,y) =
$$
$$
\sum_{l=1}^{L}\left(  \kappa_l(x) y^{2h(l)+2r(l)} + \tau_l(x) y^{h(l)+2r(l)} + \eta_l(x)y^{2r(l)} \right). 
$$
The function  $x \phi_x + y  \phi_y$ is the sum of the corresponding expressions, computed for any $l$. One has, omitting the dependence on $l$ in the last sum,
$$
x \phi_x + y  \phi_y = \sum_{l=1}^{L} x {P_l}_x + y {P_l}_y =
$$
$$
\sum \left( (x\kappa' + 2(h+r) \kappa )y^{2h+2r} + (x \tau' + (h+2r) \tau ) y^{h+2r} + (x\eta' + 2r\eta)y^{2r}  \right).
$$
For every $l$, the summands' sign is that of $(x\kappa' + 2(h+r) \kappa )y^{2h} + (x \tau' + (h+2r) \tau ) y^{h} + (x\eta' + 2r\eta) $. Performing the substitution $z=y^h$ one can study the sign of  such a $y$-trinomial by studying that of $(x\kappa' + 2(h+r) \kappa )z^2 + (x \tau' + (h+2r) \tau ) z + (x\eta' + 2r\eta) $. Its discriminant is just 
$$
(x \tau' + (h+2r) \tau ) ^2 - 4 (x\eta' + 2r\eta) (x\kappa' + 2(h+r) \kappa ).
$$
The condition $(T^+)$ implies that such a discriminant is non-positive, and that the leading term of the $y$-polynomial is non-negative, hence the trinomial is non-negative for all $y$. As a consequence, one has $x \phi_x + y  \phi_y \geq 0$ in $(a,b)$.
The set $x \phi_x + y  \phi_y = 0$ is the union of the $x$-axis together with the family of vertical  lines defined by the equalities $xf_{2k+1}'(x) + 2k f_{2k+1}(x) = 0$. Under the choosen hypothesis on the sequence $x_m$, no orbit meeting a point of $x \phi_x + y  \phi_y = 0$ can remain in such a set, since $\dot x = y$. Hence we can apply theorem \ref{teorema}.
\hfill$\clubsuit$

We can provide an example satisying the hypotheses of corollary \ref{corollarioT}. Let us set
$$
\phi(x,y)  =   (x^2 + 1)y^2+\frac{x^2}{10}y +x^2 -1  .
$$
Here there is just one $y$-trinomial, where one has $r=0$, $h=1$, $\kappa(x) = x^2+1$, $\tau(x) = \frac{x^2}{10}$, $\eta(x) = x^2-1$. One has
$$
(x\kappa' + 2(h+r) \kappa )y^{2h+2r} + (x \tau' + (h+2r) \tau ) y^{h+2r} + (x\eta' + 2r\eta)y^{2r} =
$$
$$
(4x^2 + 2 )y^2 + \frac{3x^2}{10} y  + 2 x^2 .     
$$
the discriminant of the $y$-trinomial $(4x^2 + 2 )y^2 + \frac{3x^2}{10} y  + 2 x^2$ is $-\frac{3191}{100} x^4 - 16 x^2$, which is everywhere negative, but at $0$, where it vanishes. Hence the above $y$-trinomial is positive for all $x \neq 0$, and non-negative for $x = 0$, so that the corresponding system has at most one limit cycle (see Figure 2). 

\begin{figure}[h!]
  \caption{The system $\dot x = y , \quad \dot y = -x -y( (x^2 + 1)y^2+\frac{x^2}{10}y +x^2 - 1)$ has just one limit cycle. This picture shows that the angular velocity changes along some orbits.}
  \centering
    \includegraphics[width=0.5\textwidth]{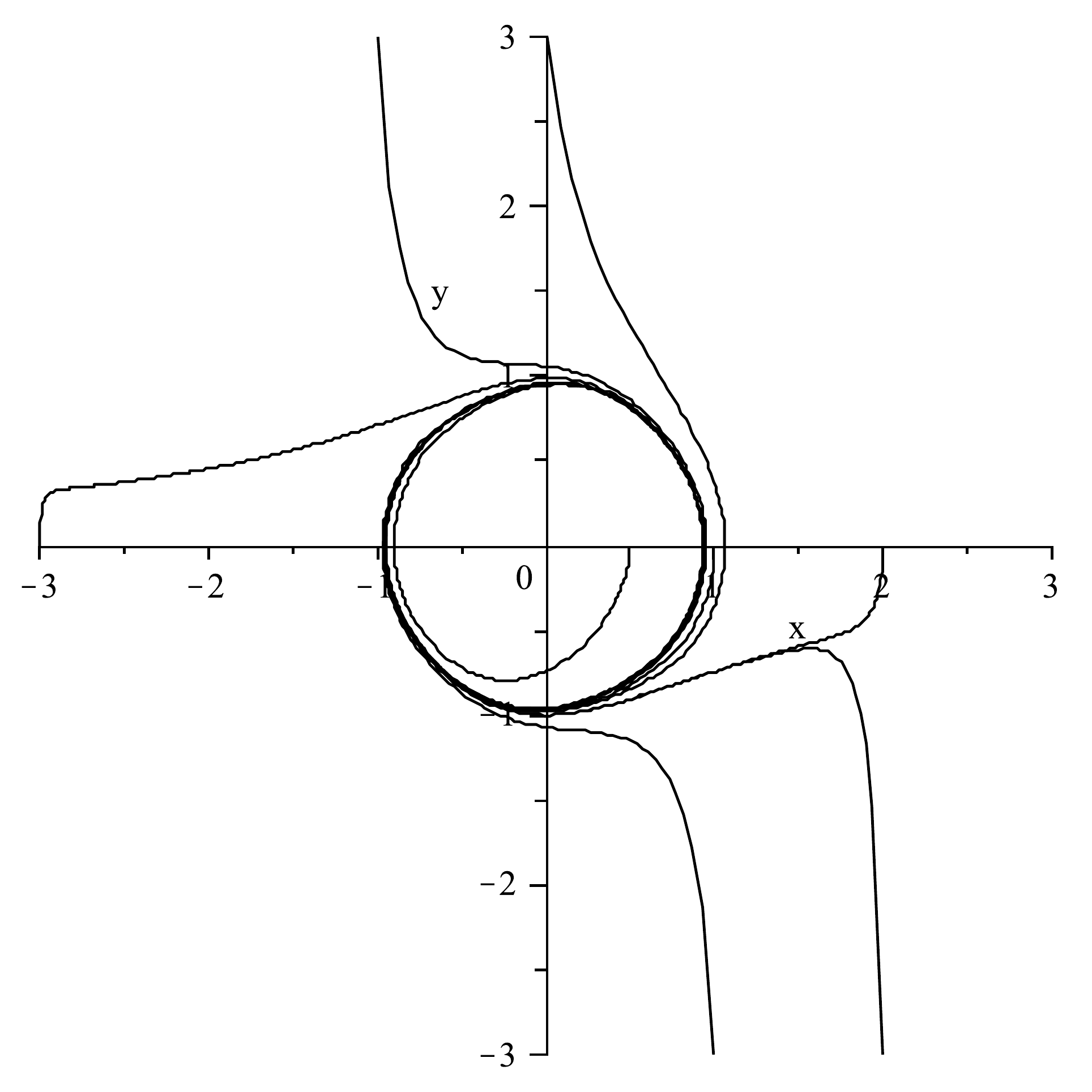}
\end{figure}

A statement similar  to corollary \ref{corollarioT} can be proved for systems satisfying the symmetric condition

${\bf (T^-)}$ for all $x \in (a,b)$, $x \neq 0$, one has $(x \tau' + (h+2r)\tau)^2 - 4( x\eta' + 2r \eta)(x \kappa' + (2h+2r) \kappa)) \leq 0$, and $x\eta' + 2r \eta \leq 0$ (or, equivalently, $x \kappa' + (2h+2r) \kappa \leq 0)$.

Next corollary is a special case of corollary  \ref{corollarioT}.

\begin{corollary}\label{corollarioCRV} Assume $f_{j}(x) \equiv 0 $ for $j$ even, $xf_{j}'(x) + (j-1) f_{j}(x) \geq 0$ for $j$ odd and for $x \in (a,b)$. If (Seq) holds, then  the system (\ref{sysCRV}) has at most one limit cycle, which is hyperbolic.
\end{corollary}
{\it Proof.} 
The function
$$
\phi(x,y) = \sum_{j=1}^{J}f_{j}(x)y^{j-1} =  \sum_{k=0}^{N}f_{2k+1}(x)y^{2k}, 
$$
can be considered as the sum of $y$-trinomials of the type $\kappa(x) y^{2h+2k} + \tau(x) y^{h+2k} + \eta(x)y^{2k}$, with $\kappa(x) =  \tau(x) \equiv 0$, $\eta(x) = f_{2k+1}(x)$.  As for the two  inequalities of $(T^+)$, one has  $(x \tau' + (h+2k)\tau)^2 - 4( x\eta' + 2k \eta)(x \kappa' + (2h+2k) \kappa)) = 0$, and $x\eta' + 2k \eta = x f_{2k+1}'(x) + 2k  f_{2k+1}(x) \geq 0$ by hypothesis. Then the conclusion comes from corollary \ref{corollarioT}.
\hfill$\clubsuit$

\begin{remark} \label{segno} The corollary \ref{corollarioCRV} is a proper extension to the uniqueness part of  theorem 1.3 in \cite{CRV}, concerned with the system (\ref{sysCRV}). 
In fact, under the hypotheses assumed in theorem 1.3 of \cite{CRV}:  

 $(L2)$ $f_{2k+1}(x) \geq 0$, for $k = 1, \dots, N$, for all $x$,  

\indent  $(L3)$ $f_{2k+1}(x)$, for $k = 0, \dots, N$, increasing for $x > 0$, decreasing for $x < 0$,  

\noindent  one has $xf_{2k+1}'(x) + 2k f_{2k+1}(x) \geq 0$, as in corollary \ref{corollarioCRV}.    \\
\indent Vice-versa, under our hypothesis $(L2)$ holds, but $(L3)$ does not necessarily hold. In order to prove $(L2)$, let us consider the half-line $x \geq 0$. First, observe that $f_{2k+1}(0) \geq 0$. Then, assume by absurd the existence of $x^*$ such that $f_{2k+1}(x^*) < 0$. Then $x^*f_{2k+1}'(x^*) \geq - 2k f_{2k+1}(x^*) > 0$. Hence $f_{2k+1}$ has a minimum at a point $x_m \in (0,x^*)$, where $f_{2k+1}'(x_m) = 0$. This contradicts $x_mf_{2k+1}'(x_m) \geq - 2k f_{2k+1}(x_m) > 0$. One works similarly on the half-line $x \leq 0$.   \\
\indent  On the other hand, our hypothesis can be satisfied even if $(L3)$ does not hold.  An example is provided by
$$
f_{2k+1}(x) = x^4 - x^2 +1, \qquad k > 0.
$$
One has $f_{2k+1}'(x) = 4x^3-2x$, so that $f_{2k+1}(x)$ is neither increasing for $x > 0$, nor decreasing for $x < 0$. Also, one has  $xf_{2k+1}'(x) + 2kf_{2k+1}(x) = (4+2k)x^4 - (2+2k)x^2 + 2k > 0 $, for $k > 0$. In fact, such a polynomial's discriminant is $\Delta = 4-24k-12k^2 < 0$ for all $k>0$, so that $xf_{2k+1}'(x) + 2kf_{2k+1}(x) $ does not vanish for any real $x$, for $k > 0$.  \\
\indent The coefficient $f_1(x)$ is not subject to the same argument, since for $k = 0$ one has $xf_{2\cdot 0+1}'(x) + 2\cdot 0 \cdot f_{2\cdot 0+1}(x) = xf_{1}'(x) \geq 0$, which does not imply any condition on the sign of $f_1(x)$. 
\end{remark}

\indent In order to prove the existence of limit cycles, we need a stronger hypothesis on the system (\ref{sysphi}).   We first prove a result about the solutions' boundedness, for systems defined on all of $\R^2$. We denote by $D_M$ the disk $\{ (x,y) : x^2+y^2 \leq M^2  \}$. We set $Z_\phi = \{ (x,y) : \phi(x,y) =0 \}$.

\begin{lemma}\label{bddness} Let $\Omega = \R^2$. If there exists  $M > 0$ such that   $\phi(x,y) \geq 0$ for all $(x,y) \not\in D_M$, and the set $Z_\phi \setminus D_M $ does not contain any non-trivial positive semi-orbit of  (\ref{sysphi}), then every solution of (\ref{sysphi}) definitely enters the disk $D_M$ and remains inside it.
\end{lemma}
{\it Proof.}
Let us consider the function $V(x,y) =   \frac 12 \left( x^2 + y^2 \right)$. Its derivative along the solutions of (\ref{sysphi}) is 
$$
\dot V(x,y) = -y^2 \phi(x,y).
$$
$\dot V(x,y) \leq 0$ out of the compact set $D_M$, hence every solution is bounded. Moreover, working as in the proof of LaSalle invariance principle (see \cite{V} for the invariance principle, \cite{S} for the details of such an argument), one can show that the positive limit set of every solution remaining in the complement of $D_M$ is contained in the set $\dot V(x,y) = 0$. Such a  positive limit set is positively invariant, i.e. if it contains a point $(x^*,y^*)$, then it contains the whole positive semi-orbit starting at $(x^*,y^*)$.  The set $\dot V(x,y) = 0$ is the union of the sets $y=0$ and $\phi=0$. By hypothesis, such a set does not contain any positive semi-orbit, hence every orbit eventually meets the set $D_M$. Since $\dot V(x,y) \leq 0$ for $x^2+y^2 \geq M^2$, every orbit definitely enters the disk $D_M$ and remains inside it.
\hfill$\clubsuit$

\begin{lemma} \label{bddness2} Let $\Omega = \R^2$. If there exists  $M > 0$ such that $\phi(x,y) \geq 0$ for all $(x,y) \not\in D_M$, and a half-line  
$ \{ (r\cos \theta^*,r \sin  \theta ^*) : r > 0\}$, such that  \\ $\phi(r\cos \theta^*,r \sin  \theta ^*) > 0 $ for  $r > M$, then every solution of (\ref{sysphi}) definitely enters the disk $D_M$ and remains inside it.
\end{lemma}
{\it Proof.} 
The derivative of $V(x,y) $ along the solutions of (\ref{sysphi}) is 
$$
\dot V(x,y) = -y^2 \phi(x,y).
$$
The solutions' boundedness comes as in  lemma \ref{bddness}. On the set $\phi(x,y) = 0$ one has $\dot x=y$, $\dot y= -x$, hence every orbit $\gamma$ starting at a point of  $\phi = 0$ is contained in a circle centered at $O$ until it meets a point where $\dot V \neq 0$. Every circle with radius greater than $M$ meets the half-line $ \{ (r\cos \theta^*,r \sin  \theta ^*)  : r > M\}$ at a point $(x_0,y_0)$. 
If $\sin  \theta ^* \neq 0$, then $\dot V(x_0,y_0) < 0$. 
If $\sin  \theta ^* = 0$, then $\phi(x,y) > 0$ in a neighbourhood of $(x_0,y_0)$, and $\dot y = x \neq 0$ implies that $\gamma$ contains a point $(x_1,y_1)$, close to  $(x_0,y_0)$, such that $\dot V(x_1,y_1) < 0$. In both cases, $\gamma$ leaves the circle pointing towards the origin. Then we may apply lemma  \ref{bddness}.
\hfill$\clubsuit$

As a particular case, we consider functions $\phi(x,y)$ obtained as sums of $y$-trinomials. Let us introduce the following definition for a $y$-trinomial $P(x,y) = \kappa(x) y^{2h+2r} + \tau(x) y^{h+2r} + \eta(x)y^{2r}$,

${\bf (T^{++})}$   there exists $\varepsilon > 0$, such that for all $x \in \R$, $|x | >  \varepsilon$,  one has $\tau(x)^2 - 4 \eta(x) \kappa(x) \leq 0$, and $ \kappa(x) \geq 0$ (or, equivalently, $\eta(x) \geq 0)$.

\begin{corollary}   \label{corbdd} Let $\Omega = \R^2$. Assume $\sum_{j=1}^{J}f_{j}(x)y^{j-1}  = \sum_{l=1}^{L} P_l(x,y) $, with $P_l(x,y) $ $y$-trinomial satisfying the condition $(T^{++})$, $l=1, \dots, L$. If, for $|x | \leq  \varepsilon$ one has $\kappa_l(x) > 0$, $l=1, \dots, L$,
then there exists a disk $D_M$ such that every solution of (\ref{sysypol}) definitely enters $D_M$ and remains inside it.
 \end{corollary}
{\it Proof.} Let us set 
$Z_l (x) = 1+ \max  \left\{ \frac{ |\tau_l(x)|}{|\kappa_l(x)|} , \frac{ |\eta_l(x)|  }{|\kappa_l(x)|}  	\right \}$. 
By Cauchy theorem about the polynomial roots, every $z$-root of $\kappa_l(x) z^2 + \tau_l(x) z + \eta_l(x)$ is contained in the $z$-interval $[-Z_l(x),Z_l(x)]$. Since $\kappa_l(x)$ is continuous and positive in $[-\varepsilon, \varepsilon]$, $Z_l(x)$ is a continuous function in $[-\varepsilon, \varepsilon]$. 
Let us set $Z_l^{max}= \max \{ Z_l(x) : x \in  [-\varepsilon,\varepsilon] \}$, and $\overline{Z} = \max  \{ Z_l^{max} : l=1, \dots , L \}$. All the zeroes of the functions $\kappa_l(x) z^2 + \tau_l(x) z + \eta_l(x)$ are contained in the rectangle 
$[-\varepsilon,\varepsilon] \times [-\overline{Z} , \overline{Z} ]$. In other words, every function $\kappa_l(x) z^2 + \tau_l(x) z + \eta_l(x)$ is positive for $x \in [-\varepsilon,\varepsilon] $ and $z\not \in  [-\overline{Z} , \overline{Z} ]$. 
As a consequence,  there exists $\overline{Y} > 0$ such that  every trinomial $P_l(x,y)$ is positive for $x \in [-\varepsilon,\varepsilon] $ and $y \not \in  [-\overline{Y},\overline{Y} ]$. 
Moreover, by $(T^{++})$, for all $x \not \in [-\varepsilon,\varepsilon]$, every $y$-trinomial is non-negative. 
Then $\phi(x,y) = \sum_{j=1}^{J}f_{j}(x)y^{j-1}  = \sum_{l=1}^{L} P_l(x,y) \geq 0 $ out of $[-\varepsilon,\varepsilon] \times [-\overline{Y} , \overline{Y} ]$.   \\
\indent  In order to check the non-positive-invariance of the set $\dot V(x,y) = 0$ it is sufficient to consider that $\phi(0,y) > 0$, for $|y| > \overline{Y}$, and apply the lemma \ref{bddness2}.
 \hfill$\clubsuit$

In next corollary we consider again the system (\ref{sysCRV}). 

\begin{corollary} \label{corCRV} Let $\Omega = \R^2$.
Assume $f_{2k}(x) = 0 $, $xf_{2k+1}'(x) + 2k f_{2k+1}(x) \geq 0$ for $x \in \R$, $k=1, \dots, N$. 
Assume additionally $ f_{1}(x) \geq  0$ for $x \in \R \setminus [-\varepsilon,\varepsilon]$, for some $\varepsilon > 0 $, and that there exists $1 \leq \overline{k} \leq N$ such that $f_{2 \overline {k}+1}(x) > 0$ on $[-\varepsilon,\varepsilon]$. 
Then there exists a disk $D_M$ such that every solution of (\ref{sysCRV}) definitely enters $D_M$ and remains inside it.
 \end{corollary}
{\it Proof.}  By remark \ref{segno}, for $k = 1, \dots, N$ one has $f_{2k+1}(x) \geq 0$, hence
 $$
 \phi(x,y) =  f_1(x) + \sum_{k=1}^N f_{2k+1}(x) y ^{2k} \geq   f_1(x) +  f_{2 \overline {k}+1}(x) y ^{2 \overline {k}}   .
 $$ 
 The function $\phi(x,y) $ is non-negative for $x\not \in (-\varepsilon,\varepsilon)$. Moreover, working as in corollary \ref{corbdd} , one can prove the existence of $\overline{Y} > 0$ such that $\phi(x,y) > 0$ for $x \in [-\varepsilon,\varepsilon] $ andl $y \not \in  [-\overline{Y},\overline{Y} ]$. Hence we may apply the lemma  \ref{bddness2}, since $\phi(0,y) > 0$ for $|y| > \overline{Y} $. 
 \hfill$\clubsuit$

\begin{remark} Our hypotheses are stronger than those of section 4.2 in \cite{CRV}, since we ask $f_1(x)$ to change sign both for positive and for negative $x$. Actually, the argument in section 4.2 of \cite{CRV} is incorrect. In fact, choosing $V(x,y) = x^2+y^2$, one has 
$$
\dot V(x,y) = -2y^2 \left( f_1(x) + \sum_{k=1}^N f_{2k+1}(x) y ^{2k} \right) .
$$
The geometrical argument in section 4.2 of  \cite{CRV} is equivalent to asking that for $r$ large enough,  $\dot V(x,y) \leq 0$. This can be false if $f_1(x)$ does not assume positive values for large $|x|$. In fact, if $f_1(x) < 0$ for large $|x|$ (see Remark 1.4 in \cite{CRV}), the curve $\dot V(x,y) = 0$ may consist of different, unbounded branches, which do not bound a positively invariant topological annulus, as needed by Poincar\'e-Bendixson theorem. This is what occurs choosing
 $$
 f_1(x) = - e^{-x^2} < 0, \qquad f_3(x) =1  - e^{-x^2} > 0.
 $$
In fact, in such a case the curve $\dot V (x,y) = 0$ consists of four unbounded connected components, separating the region $\dot V (x,y) > 0$ from the four connected regions where $\dot V (x,y) < 0$.  \\
\end{remark}

Finally, we prove a result of existence and uniqueness for the system \ref{sysphi}.

\begin{theorem}\label{teorema2} If the hypotheses of theorem \ref{teorema} and lemma \ref{bddness} hold, and $\phi(0,0) < 0$, then the system (\ref{sysphi}) has exactly one limit cycle, which is hyperbolic and  attracts every non-constant solution.
\end{theorem}
{\it Proof.} As in theorem 2 of \cite{S}
\hfill$\clubsuit$

As a consequence, we may apply the corollaries \ref{corollarioT} and  \ref{corbdd} in order to prove the limit cycle's existence and uniqueness for a special class of systems. 

\begin{corollary}  Assume   $\sum_{j=1}^{J}f_{j}(x)y^{j-1}$ to be the sum of $y$-trinomials satisfying the conditions $(T^+)$, $(T^{++})$ and (Seq).
 If $\kappa_l(x) > 0$, for $|x | \leq  \varepsilon$, $l=1, \dots, L$,
then  the system (\ref{sysphi}) has exactly one limit cycle, which is hyperbolic and attracts every non-constant solution.
\end{corollary}
{\it Proof.}  An immediate consequence of corollaries \ref{corollarioT} and  \ref{corbdd}, and theorem \ref{teorema2}.
\hfill$\clubsuit$

Finally, we have a similar result for the system (\ref{sysCRV}). 

\begin{corollary} \label{ultimo}Assume $f_{2k}(x) = 0 $, $xf_{2k+1}'(x) + 2k f_{2k+1}(x) \geq 0$ for $x \in \R$, $k=1, \dots, N$. 
Assume additionally $ f_{1}(0) < 0$, $ f_{1}(x) > 0$ for $x \in \R \setminus [-\varepsilon,\varepsilon]$, for some $\varepsilon > 0 $, and that there exists $1 \leq \overline{k} \leq N$ such that $f_{2 \overline {k}+1}(x) > 0$ on $[-\varepsilon,\varepsilon]$. 
 Then the system (\ref{sysphi}) has exactly one limit cycle, which is hyperbolic and attracts every non-constant solution.
\end{corollary}
{\it Proof.}  A straightforward consequence of theorem \ref{teorema2} and corollary \ref{corCRV}.
\hfill$\clubsuit$

\section{Non-linear $g(x)$ }

This section is similar to section 3 in \cite{S}. It contains some modifications due to the weaker hypotheses assumed in the previous section on $\phi(x,y)$. For the reader's convenience, we recall some parts of   section 3 in \cite{S}. Then we state the main result under the new, weaker  conditions,  and deduce the corresponding conditions on the system (\ref{sysphig}).

 Let us consider the equation
\begin{equation}\label{equaphig}
\ddot x + \dot x \Phi(x,\dot x) + g(x) = 0 .
\end{equation}
We assume that  $ xg(x) > 0 {\rm\ for \ } x \neq 0$, $g \in C^1((a,b),\R)$, with $a < 0 < b$, $g'(0) > 0$. We also admit the limit case $a = -\infty$ and/or $b = +\infty$. 
The main tool is the so-called Conti-Filippov transformation, which acts on the equivalent system
\begin{equation}\label{sysphig}
\dot x = y, \qquad  \dot y = - g(x) - y \Phi(x,y),
\end{equation}
in such a way to take the conservative part of the vector field into a linear one. Let us set $G(x) = \int_0^x g(s) ds$, and denote by $\sigma(x)$ the sign function, whose value is $-1$ for $x < 0$, $0$ at $0$, $1$ for $x > 0$.  Let us define the function $\alpha : \R \rightarrow \R$ as follows:
$$
 \alpha(x) =  \sigma(x) \sqrt{2G(x)}.
$$
Then Conti-Filippov transformation is the following one,
\begin{equation}\label{transf}
(u,v) = \Lambda(x,y) = (\alpha(x),y).
 \end{equation}
 $\Lambda$ transforms $(a,b)$ onto $(u^-,u^+) = \left( -\lim_{x\rightarrow a^+ } \sqrt{2G(x)},\lim_{x\rightarrow b^- }\sqrt{2G(x)} \right)$. In general, even if $(a,b) = \R$,  $\Lambda$ does not transform $\R$ onto $\R$. This occurs if and only if $\lim_{x\rightarrow -\infty } \sqrt{2G(x)} = \lim_{x\rightarrow -\infty } \sqrt{2G(x)} = +\infty$.\\  
 \indent  Since we assume $g \in C^1((a,b),\R)$, we have $\alpha  \in C^1((a,b),\R)$. The function $u = \alpha (x)$ is invertible, due to the condition $xg(x) > 0$. Let us call $x =\beta(u)$ its inverse. The condition $g'(0) > 0$ guarantees the differentiability of $\beta(u)$ at $O$. For $x \neq 0$, or, equivalently, for $u \neq 0$, one has,
\begin{equation}\label{derivate}
\alpha'(x) =  \frac{\sigma(x) g(x)}{\sqrt{2G(x)}}, \qquad \beta'(u) =  \frac{1}{\alpha'(\beta(u))} =  \frac{\sigma(x)\sqrt{2G(x)}}{g(x)} = 	\frac{u}{g(\beta(u))}.
 \end{equation}
 For $x=u=0$ one has, 
$$
 \alpha'(0) = \sqrt{g'(0)}, \qquad  \beta'(0) =  \sqrt{\frac{1}{g'(0)}} .
$$
Finally,
\begin{equation}\label{root}
\lim_{u \rightarrow 0} \frac{g(\beta(u))}{u} =  \sqrt{g'(0)} > 0.
\end{equation}

In next theorem we consider only the condition corresponding to $x\phi_x + y\phi_y \geq 0$. The alternative case, $x\phi_x + y\phi_y \leq 0$, can be treated similarly. Let us set

$$
\Psi =  \frac {\sigma \sqrt{2G} }{g} \left[ 2G \frac{ \Phi_x g - \Phi g' }{g ^2} +  \Phi   \right]
+ y \Phi_y.
$$
\bigskip
\begin{theorem}\label{teorema-n} Assume $g \in C^1((a,b),\R)$, with $g'(0) > 0$ and $xg(x) > 0$ for $x\neq 0$. If  $\Psi(x,y) \geq 0$ and the set $\Psi(x,y) = 0$ does not contain any non-trivial orbit, then (\ref{sysphig}) has at most one limit cycle in the region $(a,b) \times \R$, which is hyperbolic.
\end{theorem}
{\it Proof.} As in the proof of theorem 3 in \cite{S}, one transforms the system (\ref{sysphig}) into a system of the form (\ref{sysphi}), by means of the Conti-Filippov transformation.  \\
\indent  For $u \neq 0$, the transformed system has the form
\begin{equation}\label{1}
\dot u = v \frac{g(\beta(u))}{u}, \qquad  \dot v = -  g(\beta(u)) -  v \Phi(\beta(u),v) ,
\end{equation}
We may multiply the system (\ref{1}) by $ \frac{u}{g(\beta(u))} $, obtaining a new system having the same orbits as (\ref{1}), 
\begin{equation} \label{CFphi}
\dot u = v ", \qquad  \dot v = -  u -  v  \frac{u \Phi(\beta(u),v)}{g(\beta(u))} =  -  u -  v \phi(u,v),
\end{equation}
By (\ref{root}), the new system is regular also for $x=0$. As proved in \cite{S}, one has $\Psi(\beta(u),v) = u\phi_u(u,v) + v\phi_v(u,v) $, hence the hypothesis  $\Psi(x,y) \geq 0$ implies the star-shapedness condition on $\phi(u,v)$. Moreover, the set $\Psi(x,y) = 0$ is transformed into the set $u\phi_u + v\phi_v = 0$.
Then one applies theorem \ref{teorema} to complete the proof.
\hfill$\clubsuit$

In order to apply the above theorem to $y$-polynomial equations, as in the previous section's corollaries, it is convenient to write the $y$-coefficients form, after the application of Conti-Filippov transformation. If 
$$
\Phi(x,y) = \sum_{j=1}^{J}f_{j}(x)y^{j-1}
$$
then 
$$
\phi(u,v) = \frac{u \Phi(\beta(u),v)}{g(\beta(u))} =  \sum_{j=1}^{J}  \frac {u f_{j}(\beta(u))} {g(\beta(u))} v^{j-1} = 
 \sum_{j=1}^{J}\tilde{  f}_{j}(u) v^{j-1} ,
$$
where $ \tilde{  f}_{j}(u)  =  \frac {u f_{j}(\beta(u))} {g(\beta(u))}$.
If $\Phi(x,y)$ is the sum of $y$-trinomials $\kappa_l(x) y^{2h+2r} + \tau_l(x) y^{h+2r} + \eta_l(x)y^{2r}$, then $\phi(u,v)$ is the sum of $v$-trinomials, whose coefficients $\tilde{\kappa_l}(u) v^{2h+2r} + \tilde{\tau_l}(u) v^{h+2r} + \tilde{\eta_l}(u)v^{2r}$are obtained from the original ones in a similar way,
$$
\tilde{\kappa_l}(u ) =  \frac {u \kappa_l(\beta(u))} {g(\beta(u))} ,  \qquad
\tilde{\tau_l}(u ) =  \frac {u \tau_l(\beta(u))} {g(\beta(u))} , \qquad
\tilde{\eta_l}(u ) =  \frac {u \eta_l(\beta(u))} {g(\beta(u))} .
$$
Writing the condition $(T^+)$ for such polynomials generates quite cumbersome expressions. In next corollary we only write the conditions for the simplest case, that of an odd $y$-polinomial.   \\
\indent Given a function $f: I \rightarrow \R$, consider the condition     

${\bf (H_g^{j})}$ \quad  $\forall x \neq 0$: \  $x \left[  j f(x)g(x) + 2G(x) \left(  \frac{f'(x)g(x) - f(x)g(x)'}{g(x)} \right) \right] \geq 0. $

If the above inequality holds strictly for $x \neq 0$, i.e. if the left-hand side does not vanish for $x\neq 0$, we say that the strict condition ($H_g^{j}$) holds at $x$.

\begin{corollary}\label{corCRV-n} Assume $g \in C^1((a,b),\R)$, with $g(x)= x + o(x)$ and $xg(x) > 0$ for  $x\neq 0$. Let $f_{2k}\equiv 0$ and $f_{2k+1}(x)$ satisfy the condition ($H_g^{2k+1}$) for $k=0, \dots, N$. If there exists a sequence $x_m$ converging to $0$, such that for every  $m$ there exists $k(m)$ satisying the strict ($H_g^{2k+1}$) condition at $x_m$, then (\ref{sysypol}) has at most one limit cycle in the region $(a,b) \times \R$, which is hyperbolic.
\end{corollary}
{\it Proof.}
 As observed above, the transformed system has the form \begin{equation}\label{sysurep}
\dot u = v  , \qquad  \dot v = - u -  v\sum_{k=0}^{N}  \frac{ u f_{2k+1}(\beta(u))}{g(\beta(u))} v^{2k} .
\end{equation}
Applying corollary \ref{corollarioCRV} to such a system requires to check, for $k = 0, \dots, N$, the condition
$$
u\left(  \frac{ u f_{2k+1}(\beta(u))}{g(\beta(u))} \right)' + 2k  \frac{ u f_{2k+1}(\beta(u))}{g(\beta(u))}  \geq 0.
$$
Performing standard computations and recalling that $u =   \sigma(x) \sqrt{2G(x)}$, one proves that, for $u\neq 0$, the inequality
$$
u\left(  \frac{ u f(\beta(u))}{g(\beta(u))} \right)' + 2k  \frac{ u f(\beta(u))}{g(\beta(u))}  \geq 0
$$
is equivalent to 
$$
\sigma(x) \sqrt{G(x)} \left[  (2k+1) f(x)g(x) + 2G(x) \left(  \frac{f'(x)g(x) - f(x)g(x)'}{g(x)} \right) \right]  \geq 0,
$$
for $x \neq 0$.
The sign of $\sigma(x) \sqrt{G(x)}$ is the same as that of $x$, hence the last inequality is equivalent to  ($H_g^{2k+1}$). Then one applies corollary \ref{corollarioCRV} to complete the proof.
\hfill$\clubsuit$

Now we prove the analogous of lemma \ref{bddness}. Let us set
$$
E(x,y) = G(x) + \frac {y^2}{2}
$$
For $r > 0$, we set $\Delta_r = \{ (x,y) : 2 E(x,y) < r^2 \}$. As in the previous section, we set $Z_\Phi = \{ (x,y) : \Phi(x,y) =0 \}$.

\begin{lemma}\label{bddnessg} Assume $g \in C^1((a,b),\R)$, with $\int_0^a g(x) dx = \int_0^b g(x) dx = +\infty$. 
If there exists  $M > 0$ such that   $\Phi(x,y) \geq 0$ for all $(x,y) \not\in \Delta_M$, and the set $Z_\Phi \setminus \Delta_M $ does not contain any non-trivial positive semi-orbit of  (\ref{sysphig}), then every solution of (\ref{sysphi}) definitely enters the set $\Delta_M$ and remains inside it.
\end{lemma}

{\it Proof.} Under the integral conditions in the hypotheses, the Conti-Filippov transformation is a global diffeomorpism of the region $(a,b) \times \R$ onto the plane, because
$$
\lim_{x\rightarrow a^+ } \sqrt{2G(x)} = \lim_{x\rightarrow b^- }\sqrt{2G(x)} = +\infty.
$$
The level sets of $E(x,y)$ are compact subsets of $[a,b] \times \R$, i. e. they have positive distance from its boundary. 
The sets $ \Delta_r$ are taken into the sets $D_r$.  
Then one can apply the lemma  \ref{bddness} to the system (\ref{CFphi}). In fact, the derivative of the Liapunov function $V(u,v) = \frac 12 (u^2+v^2)$ along the solutions of  (\ref{sysurep}) is just 
$$ 
\dot V (u,v) =  - v^2  \frac{u\Phi(\beta(u),v)} {g(\beta(u))} .
$$
The function $ \frac{u} {g(\beta(u))}$ is positive for $u \neq 0$, hence the hypotheses of lemma \ref{bddness} are satisfied by the system (\ref{sysurep}). As a consequence, the conclusions of lemma  \ref{bddness} hold for the system (\ref{sysurep}), and applying the inverse transformation $\Lambda^{-1}$ one obtains the thesis.
\hfill$\clubsuit$

\begin{theorem}\label{teorema2-n}  Let $g\in C^1((a,b),\R)$. If the hypotheses of theorem \ref{teorema-n} and lemma \ref{bddnessg} hold on $\R$, and $\Phi(0,0) < 0$, then the system (\ref{sysphig}) has exactly one limit cycle in $(a,b) \times \R$, which is hyperbolic and attracts every non-constant solution.
\end{theorem}
{\it Proof.} As the proof of theorem \ref{teorema2}, replacing the Liapunov function $V(x,y)$ with the Liapunov function $E(x,y)$. 
\hfill$\clubsuit$

Finally, we prove a result analogous to corollary \ref{ultimo}.

\begin{corollary} Assume $g \in C^1(\R,\R)$, $g'(0) > 0$, $xg(x) > 0$ for  $x\neq 0$, $\int_0^{\pm \infty} g(x) dx = +\infty$. Let $f_{2k}\equiv 0$ and $f_{2k+1}(x)$ satisfy the condition ($H_g^{2k+1}$) for $k=0, \dots, N$. Assume there exists a sequence $x_m$ converging to $0$, such that for every  $m$ there exists $k(m)$ satisying the strict ($H_g^{2k+1}$) condition at $x_m$.
Assume additionally $ f_{1}(0) < 0$, $ f_{1}(x) > 0$ for $x \in \R \setminus [-\varepsilon,\varepsilon]$, for some $\varepsilon > 0 $, and that there exists $1 \leq \overline{k} \leq N$ such that $f_{2 \overline {k}+1}(x) > 0$ on $[-\varepsilon,\varepsilon]$. 
 Then the system (\ref{sysypol}) has exactly one limit cycle, which attracts every non-constant solution.
\end{corollary}
{\it Proof.}  
The sign of $f_j(x)$ is the same as that of  $\tilde{f}_j(x)$, since  $\frac {u} {\beta(u)} >0 $ for $u \neq 0$. 
Then the statement is a straightforward consequence of theorem \ref{teorema2-n}  and corollary \ref{corCRV-n}.
\hfill$\clubsuit$

{\bf  \large Acknowledgements}

This paper has been partially supported by the GNAMPA 2009 project \lq\lq Studio delle traiettorie di equazioni differenziali ordinarie\rq\rq.

\end{document}